\documentclass[12pt,twoside]{article}
\usepackage[english]{babel}
\usepackage[latin1]{inputenc}
\usepackage{amsmath}
\usepackage{amssymb,amsfonts}
\usepackage{graphicx}
\usepackage{times}

\usepackage{float}

\usepackage{amsmath}
\usepackage{amssymb}
\usepackage{amscd}
\usepackage{amsthm}

\usepackage{graphicx}
\usepackage{epstopdf}

\usepackage{longtable}

\usepackage{array}
\usepackage{times,amssymb,amscd}
\newtheorem{thm}{Theorem}[section]
\newtheorem{cor}[thm]{Corollary}

\newtheorem{lem}[thm]{Lemma}

\newtheorem{defn}[thm]{Definition}
\newtheorem{rem}[thm]{Remark}
\numberwithin{equation}{section}

\newcommand{\bA}{\mathbf{A}}

\newcommand{\bE}{\mathbf{E}}

\newcommand{\bH}{\mathbf{H}}

\newcommand{\bL}{\mathbf{L}}

\newcommand{\bR}{\mathbf{R}}
\newcommand{\bS}{\mathbf{S}}
\newcommand{\bV}{\mathbf{V}}

\newcommand{\be}{\mathbf{e}}

\newcommand{\bt}{\mathbf{t}}

\newcommand{\cP}{\mathcal{P}}
\newcommand{\cS}{\mathcal{S}}
\newcommand{\cT}{\mathcal{T}}
\newcommand{\cC}{\mathcal{C}}

\newcommand{\cB}{\mathcal{B}}

\newcommand{\EUC}{\mathbf E^3}

\newcommand{\HYP}{\bH^3}
\newcommand{\SXR}{\bS^2\!\times\!\bR}
\newcommand{\HXR}{\bH^2\!\times\!\bR}
\newcommand{\SLR}{\widetilde{\bS\bL_2\bR}}
\newcommand{\NIL}{\mathbf{Nil}}
\newcommand{\SOL}{\mathbf{Sol}}

\newcommand{\Rr}{R_\Gamma^C}
\def\NN{\mathbb{N}}

\begin{document}
\pagestyle{myheadings}
\markboth{\centerline{Judit Sajtos and Jen\H o Szirmai}}
{Lattice-like Packings and Coverings $\dots$}
\title
{Lattice-like Packings and Coverings with Congruent Translation Balls and Cylinders in $\SOL$ geometry
}

\author{Judit Sajtos and Jen\H o Szirmai \\
\normalsize Department of Algebra and Geometry, Institute of Mathematics,\\
\normalsize Budapest University of Technology and Economics, \\
\normalsize M\H uegyetem rkp. 3., H-1111 Budapest, Hungary \\
\normalsize szirmai@math.bme.hu
\date{\normalsize{\today}}}

\maketitle
\begin{abstract}

The aim of this paper is to study lattice-like coverings with congruent translation balls and the packings and coverings with a type of translation 
cylinders in $\SOL$ space related to the fundamental lattices. We introduce the notions of the densities of the considered problems
and give upper estimate to ball coverings using the radii and the volumes of the circumscribed translation spheres of given {\it translation tetrahedra}.
Moreover we determine the exact optimal packing and covering densities of a type of cylinder packings belonging to the fundamental lattices. 
\end{abstract}

\section{Introduction} \label{section1}
The basic problems of packings and coverings are classical questions of geometry and they deeply connected 
to number theory, crystallography, and many other important areas of mathematics. 
Therefore, determining the densest packing and the thinnest covering with congruent copies of a given body 
have always been considered important problems. Our work is related to this topic. 

The beginning of this topic is the famous Kepler's conjecture, which concerned the densest packing of congruent balls in the Euclidean space $\EUC$. 
One major recent development has been the settling of the long-standing Kepler conjecture, 
part of Hilbert's 18th problem, by Thomas Hales at the turn of the 21st century.
Hales' computer-assisted proof was largely based on a program set forth by L. Fejes-T\'oth in the 1950's.

The second author extended the classic Kepler's problem to non-constant curvature 
Thurston geometries $\SXR,~\HXR,~$ $\SLR,~\NIL,~\SOL$, 
in  \cite{MSz}, \cite{Sz07}, \cite{Sz13-2}, \cite{Sz13-3}, \cite{Sz22-3}; from which in this paper we consider the $\SOL$ geometry. 

At present, we consider the lattice-like covering problem with congruent translation balls in $\SOL$ space.
These questions are related to the theory of the Dirichlet-Voronoi cells (briefly $D-V$ cells).
In $3$-dimensional spaces of constant curvature the $D-V$ cells are widely investigated, but in the other Thurston geometries, $\SXR$, $\HXR$,
$\NIL$, $\SOL$, $\SLR$, there are few results in this topic. Let $X$ be one of the above five geometries and $\Gamma$ be one of its discrete isometry groups.
Moreover, we distinguish two distance function types: $d^g$ is the usual geodesic distance function and $d^t$ is the translation distance function (see Section 2).
Therefore, we obtain two types of the $D-V$ cells regarding the two distance functions. {\it In this paper we use the notion of the translation distance.} 

The investigation of this issue brought many interesting results and opened an important 
path in the direction of non-Euclidean crystal geometry (see the survey \cite{Sz22-3} 
and \cite{G--K--K}, \cite{Fe23}, \cite{FTL}, \cite{MSzV}, \cite{MSz14}, \cite{MSzV17}, \cite{Sz23}, \cite{Sz07}, \cite{Sz12-1}, \cite{Sz13-1}, \cite{Sz13-2}, 
\cite{Sz25}, \cite{T}, \cite{S}, \cite{YSz24-1}, \cite{YSz24-2}).
We mention only some here:								        
\begin{enumerate}
\item A candidate of the densest geodesic ball
packing is described in \cite{Sz13-2}. In the Thurston geometries the greatest known density was $\approx 0.8533$,
which is not realized by a packing with {\it equal balls} of the hyperbolic
space $\HYP$. However, it is attained, e.g., by a {\it horoball packing} of
$\overline{\bH}^3$ where the ideal centres of horoballs lie on the
absolute figure of $\overline{\bH}^3$ inducing the regular ideal
simplex tiling $(3,3,6)$ by its Coxeter-Schl\"afli symbol.
In \cite{Sz13-1} we have presented a geodesic ball packing in the $\SXR$ geometry
whose density is $\approx 0.8776$.
\item
In \cite{MSz12} we classified $\SOL$ lattices in an algorithmic way into 17 (seventeen) types, 
in analogy with the 14 Bravais-types of the Euclidean 3-lattices, but infinitely many $\SOL$ affine 
equivalence classes, in each type.
This way the discrete isometry groups of compact fundamental domains (crystallographic groups) can also be classified into infinitely 
many classes but finitely many types, left to other publication.
For this we shall study relations between 
$\SOL$ lattices and lattices of the pseudoeuclidean (or here rather called Minkowskian) plane \cite{AQ}. 

Moreover, the densest translation ball packing by so-called fundamental lattices, which is one (Type {\bf I/1}) of the 17 Bravais-type
of $\SOL$-lattices described in \cite{MSz06} is determined in \cite{Sz13-3}. It turns out that the optimal arrangement has a richer symmetry group (in Type {\bf I/2}) 
for $N=4$. This density is $\delta \approx 0.56405083$ and the kissing number of the balls
to this packing is 6.
\item 
Another important direction is the study of Apollonius surfaces, a special case of 
which are the bisector surfaces. These are of fundamental importance in the determinations of Dirichlet-Voronoi cells.
The geodesic-like Apollonius surfaces are investigated in $\SXR$, $\HXR$ and $\NIL$ geometries in \cite{Sz22, Sz23-1}, and the translation-like Apollonius surfaces in Thurston geometries in \cite{Cs-Sz25}.
Moreover, we studied the translation-like equidistant surfaces in $\SOL$ and $\NIL$ geometries in \cite{Sz17-1} and in \cite{SzV19}.
In \cite{PSSz10}, \cite{PSSz11}, \cite{PSSz11-1} we studied the geodesic-like bisector surfaces in $\SXR$ and $\HXR$ spaces.
\end{enumerate}
Among Thurston geometries with non-constant curvature, even fewer results are known regarding 
coverings by congruent translation or geodesic balls:
\begin{enumerate}
\item
In \cite{Sz13-1} the second author investigeted {\it lattice-like ball coverings} in $\NIL$ space, introduced the notion of the density of the considered coverings and
gave upper and lower estimates to it, moreover formulated a conjecture for the ball arrangement of the least dense lattice-like geodesic
ball covering and gave its covering density $\Delta\approx 1.42900615$. 
\item
In \cite{SzV19} A. Vr\'anics and the second author studied the $\NIL$ translation ball coverings, determined the equation of the translation-like bisector surface
of any two points, proved that the isosceles property of a translation triangle is not equivalent to two angles of the triangle being equal and that
the triangle inequalities do not remain valid for translation triangles in general.
We developed a method to determine the centre and the radius of the circumscribed translation sphere of a given {\it translation tetrahedron}.

Moreover, we investigated lattice-like coverings with congruent translation balls in $\NIL$
space, introduced the notion of the density of the considered coverings
and gave upper estimate to it using the radius and the volume of the circumscribed translation sphere of a given {\it translation tetrahedron}. 
The found minimal upper bound density of the translation ball coverings is $\Delta \approx 1.42783$.
\end{enumerate}

{\it In this work, related to this coverage topic,} we develop an algorithm to determine the least dense ball covering arrangement of a given periodic parallelepiped 
$\SOL$ tiling using the fact that this fundamental domain of the tiling can be decomposed into tetrahedra.
We will use the method described in the paper \cite{Sz17-1} to determine the centre and the radius of the circumscribed translation sphere of a given {\it translation tetrahedron}.
Applying the above procedure we determine the minimal covering density of some lattice types and thus we give an upper bound of the lattice-like covering 
density related to the important, so called fundamental lattice types (and thus for all lattice types) (see Theorems 3.4).

Moreover, we introduce a type of $\SOL$ translation cylinders, 
define their volume and introduce the notion of cylinder covering related to the lattice-like tilings and define their packing and covering densities. 
Furthermore, we determine exactly the density of their densest packings or thinnest coverings related to fundamental lattices.  
Also, for the general lattice types, we give an exact lower bound of the density of cylinder coverings and an exact upper bound of the density of the packings,
which are the same as the optimal covering and packing densities with congruent circles in the Euclidean plane (see Theorems 4.12 and 4.14).
\section{On Sol geometry}
\label{sec:1}

In this Section we summarize the significant notions and notations of real $\SOL$ geometry (see \cite{M97}, \cite{S}).

$\SOL$ is defined as a 3-dimensional real Lie group with multiplication
\begin{equation}
     \begin{gathered}
(a,b,c)(x,y,z)=(x + a e^{-z},y + b e^z ,z + c).
     \end{gathered} \tag{2.1}
     \end{equation}
We note that the conjugation by $(x,y,z)$ leaves invariant the plane $(a,b,c)$ with fixed $c$:
\begin{equation}
     \begin{gathered}
(x,y,z)^{-1}(a,b,c)(x,y,z)=(x(1-e^{-c})+a e^{-z},y(1-e^c)+b e^z ,c).
     \end{gathered} \tag{2.2}
     \end{equation}
Moreover, for $c=0$, the action of $(x,y,z)$ is determined only by its $z$-component, where $(x,y,z)^{-1}=(-x e^{z}, -y e^{-z} ,-z)$. Thus the $(a,b,0)$ plane is distinguished as a {\it base plane} in
$\SOL$, or with other words, $(x,y,0)$ is a normal subgroup of $\SOL$.
$\SOL$ multiplication can also be interpreted affinely (projectively) by "right translations"
on its points as the following matrix formula shows, according to (2.1):
     \begin{equation}
     \begin{gathered}
     (1,a,b,c) \to (1,a,b,c)
     \begin{pmatrix}
         1&x&y&z \\
         0&e^{-z}&0&0 \\
         0&0&e^z&0 \\
         0&0&0&1 \\
       \end{pmatrix}
       =(1,x + a e^{-z},y + b e^z ,z + c)
       \end{gathered} \tag{2.3}
     \end{equation}
by row-column multiplication.
This defines "translations" $\mathbf{L}(\mathbf{R})= \{(x,y,z): x,~y,~z\in \mathbf{R} \}$
on the points of space $\SOL= \{(a,b,c):a,~b,~c \in \mathbf{R}\}$.
These translations are not commutative, in general.
Here we can consider $\mathbf{L}$ as a projective collineation group with right actions in homogeneous
coordinates as usual in classical affine-projective geometry.
We will use the Cartesian homogeneous coordinate simplex $E_0(\be_0)$, $E_1^{\infty}(\be_1)$, \ $E_2^{\infty}(\be_2)$, \ 
$E_3^{\infty}(\be_3), \ (\{\be_i\}\subset \bV^4$ \ $\text{with the unit point}$ $E(\be = \be_0 + \be_1 + \be_2 + \be_3 ))$
which is distinguished by an origin $E_0$ and by the ideal points of coordinate axes, respectively.
Thus {$\SOL$} can be visualized in the affine 3-space $\bA^3$
(so in Euclidean space $\bE^3$) as well.

In this affine-projective context in \cite{M97} E. Moln\'ar has derived the usual infinitesimal arc-length square at any point
of $\SOL$, by pull back translation, as follows
\begin{equation}
   \begin{gathered}
      (ds)^2:=e^{2z}(dx)^2 +e^{-2z}(dy)^2 +(dz)^2.
       \end{gathered} \tag{2.4}
     \end{equation}
Hence we get an infinitesimal Riemann metric invariant under translations, by the symmetric metric tensor field $g$ on $\SOL$ by components as usual.

It will be important for us that the full isometry group Isom$(\SOL)$ has eight components, since the stabilizer of the origin
is isomorphic to the dihedral group $\mathbf{D_4}$, generated by two involutive (involutory) transformations, preserving (2.4):
\begin{equation}
   \begin{gathered}
      (1)  \ \ y \leftrightarrow -y; \ \ (2)  \ x \leftrightarrow y; \ \ z \leftrightarrow -z; \ \ \text{i.e. first by $3\times 3$ matrices}:\\
     (1) \ \begin{pmatrix}
               1&0&0 \\
               0&-1&0 \\
               0&0&1 \\
     \end{pmatrix}; \ \ \
     (2) \ \begin{pmatrix}
               0&1&0 \\
               1&0&0 \\
               0&0&-1 \\
     \end{pmatrix}; \\
     \end{gathered} \tag{2.5}
     \end{equation}
     with its product, generating a cyclic group $\mathbf{C_4}$ of order 4
     \begin{equation}
     \begin{gathered}
     \begin{pmatrix}
                    0&1&0 \\
                    -1&0&0 \\
                    0&0&-1 \\
     \end{pmatrix};\ \
     \begin{pmatrix}
               -1&0&0 \\
               0&-1&0 \\
               0&0&1 \\
     \end{pmatrix}; \ \
     \begin{pmatrix}
               0&-1&0 \\
               1&0&0 \\
               0&0&-1 \\
     \end{pmatrix};\ \
     \mathbf{Id}=\begin{pmatrix}
               1&0&0 \\
               0&1&0 \\
               0&0&1 \\
     \end{pmatrix}.
     \end{gathered} \notag
     \end{equation}
     Or we write by collineations fixing the origin $O=(1,0,0,0)$:
\begin{equation}
(1) \ \begin{pmatrix}
         1&0&0&0 \\
         0&1&0&0 \\
         0&0&-1&0 \\
         0&0&0&1 \\
       \end{pmatrix}, \ \
(2) \ \begin{pmatrix}
         1&0&0&0 \\
         0&0&1&0 \\
         0&1&0&0 \\
         0&0&0&-1 \\
       \end{pmatrix} \ \ \text{of form (2.3)}. \tag{2.6}
\end{equation}
A general isometry of $\SOL$ to the origin $O$ is defined by a product $\gamma_O \tau_X$, with first $\gamma_O$ of form (2.6) then $\tau_X$ of (2.3). For
a general point $A=(1,a,b,c)$, this will be a product $\tau_A^{-1} \gamma_O \tau_X$, mapping $A$ into $X=(1,x,y,z)$.

Conjugation of translation $\tau$ by an above isometry $\gamma$, as $\tau^{\gamma}=\gamma^{-1}\tau\gamma$ also denotes it, will also be used by
(2.3) and (2.6) or also by coordinates with above conventions.

We remark only that the role of $x$ and $y$ can be exchanged throughout the paper, but this leads to the mirror interpretation of $\SOL$.
As formula (2.4) fixes the metric of $\SOL$, the change above is not an isometry of a fixed $\SOL$ interpretation. Other conventions are also accepted
and used in the literature.

{\it $\SOL$ is an affine metric space (affine-projective in the sense of the unified formulation of \cite{M97}). Therefore, its linear, affine, unimodular,
etc. transformations are defined as those of the embedding affine space.}
\subsection{Translation curves}

We consider a $\SOL$ curve $(1,x(t), y(t), z(t) )$ with a given starting tangent vector at the origin $O=(1,0,0,0)$
\begin{equation}
   \begin{gathered}
      u=\dot{x}(0),\ v=\dot{y}(0), \ w=\dot{z}(0).
       \end{gathered} \tag{2.7}
     \end{equation}
For a translation curve let its tangent vector at the point $(1,x(t), y(t), z(t) )$ be defined by the matrix (2.3)
with the following equation:
\begin{equation}
     \begin{gathered}
     (0,u,v,w)
     \begin{pmatrix}
         1&x(t)&y(t)&z(t) \\
         0&e^{-z(t)}&0&0 \\
         0&0&e^{z(t)}& 0 \\
         0&0&0&1 \\
       \end{pmatrix}
       =(0,\dot{x}(t),\dot{y}(t),\dot{z}(t)).
       \end{gathered} \tag{2.8}
     \end{equation}
Thus, {\it translation curves} in $\SOL$ geometry (see \cite{MoSzi10} and \cite{MSz}) are defined by the first order differential equation system
$\dot{x}(t)=u e^{-z(t)}, \ \dot{y}(t)=v e^{z(t)},  \ \dot{z}(t)=w,$ whose solution is the following:
\begin{equation}
   \begin{gathered}
     x(t)=-\frac{u}{w} (e^{-wt}-1), \ y(t)=\frac{v}{w} (e^{wt}-1),  \ z(t)=wt, \ \mathrm{if} \ w \ne 0 \ \mathrm{and} \\
     x(t)=u t, \ y(t)=v t,  \ z(t)=z(0)=0 \ \ \mathrm{if} \ w =0.
       \end{gathered} \tag{2.9}
\end{equation}
We assume that the starting point of a translation curve is the origin, because we can transform a curve into an
arbitrary starting point by translation (2.3), moreover, unit velocity translation can be assumed :
\begin{equation}
\begin{gathered}
        x(0)=y(0)=z(0)=0; \\ \ u=\dot{x}(0)=\cos{\theta} \cos{\phi}, \ \ v=\dot{y}(0)=\cos{\theta} \sin{\phi}, \ \ w=\dot{z}(0)=\sin{\theta}; \\
        - \pi < \phi \leq \pi, \ -\frac{\pi}{2} \leq \theta \leq \frac{\pi}{2}. \tag{2.10}
\end{gathered}
\end{equation}
\begin{defn}
The translation distance $d^t(P_1,P_2)$ between the points $P_1$ and $P_2$ is defined by the arc length of the above translation curve
from $P_1$ to $P_2$.
\end{defn}
Thus we obtain the parametric equation of the the {\it translation curve segment} $t(\phi,\theta,t)$ with starting point at the origin in direction
\begin{equation}
\bt(\phi, \theta)=(\cos{\theta} \cos{\phi}, \cos{\theta} \sin{\phi}, \sin{\theta}) \tag{2.11}
\end{equation}
where $t \in [0,r], ~ r \in \bR^+$. If $\theta \ne 0$, then the system of equations is:
\begin{equation}
\begin{gathered}
        \left\{ \begin{array}{ll}
        x(\phi,\theta,t)=-\cot{\theta} \cos{\phi} (e^{-t \sin{\theta}}-1), \\
        y(\phi,\theta,t)=\cot{\theta} \sin{\phi} (e^{t \sin{\theta}}-1), \\
        z(\phi,\theta,t)=t \sin{\theta}.
        \end{array} \right. \\
        \text{If $\theta=0$ then}: ~  x(t)=t\cos{\phi} , \ y(t)=t \sin{\phi},  \ z(t)=0.
        \tag{2.12}
\end{gathered}
\end{equation}
\begin{defn} The sphere of radius $r >0$ with centre at the origin (denoted by $S^t_O(r)$) with the usual longitude and altitude parameters
$- \pi < \phi \leq \pi$,  $-\frac{\pi}{2} \leq \theta \leq \frac{\pi}{2}$, respectively, by (2.10), is specified by the equations (2.12), where $t=r$.
\end{defn}
\begin{defn}
 The body of the translation sphere of centre $O$ and of radius $r$ in the $\SOL$ space is called a translation ball, denoted by $B^t_{O}(r)$,
 i.e. $Q \in B^t_{O}(r)$ iff $0 \leq d^t(O,Q) \leq r$.
\end{defn}

In \cite{Sz13-3} we proved the volume formula of the translation ball $B^t_{O}(r)$ of radius $r$:
\begin{thm}[\cite{Sz13-3}]
\begin{equation}
\begin{gathered}
Vol(B^t_{O}(r))=\int_{V} \mathrm{d}x ~ \mathrm{d}y ~ \mathrm{d}z = \\ = \int_{0}^{r} \int_{-\frac{\pi}{2}}^{\frac{\pi}{2}} \int_{-\pi}^{\pi} \frac{\cos{\theta}}{\sin^2{\theta}}
(e^{\rho \sin{\theta}}+e^{-\rho \sin{\theta}}-2) \ \mathrm{d}\phi \ \mathrm{d}\theta \ \mathrm{d}\rho =  \\
= 4 \pi \int_{0}^{r} \int_{-\frac{\pi}{2}}^{\frac{\pi}{2}} \frac{\cos{\theta}}{\sin^2{\theta}}
(\cosh(\rho \sin{\theta})-1) \ \mathrm{d}\theta \ \mathrm{d}\rho. \notag
\end{gathered}
\end{equation}
\end{thm}
An easy power series expansion with substitution $\rho \sin{\theta}=:z$
can also be applied, no more detailed.
From the equation of the translation spheres $S^t_O(r)$ (see (2.12)) it follows that the plane sections
of these spheres, given by parameters $\theta$ and $r$, parallel to the $[x,y]$ plane are ellipses by the equations
(see Fig.~4):
\begin{equation}
\begin{gathered}
\frac{x^2}{k_1^2}+\frac{y^2}{k_2^2}=1 \ \mathrm{where} \\
k_1^2=(-\cot{\theta} (e^{-r \sin{\theta}}-1))^2, \ \ \ k_2^2=(\cot{\theta} (e^{r \sin{\theta}}-1))^2. \tag{2.13}
\end{gathered}
\end{equation}
\subsection{Translation-like bisector surfaces, translation tetrahedra and their circumscribed spheres}
The process of determining of the translational bisector surfaces was carried out in our article \cite{Sz17-1}. 
Here we recall only the results that are necessary to determine the center and radius 
of the circumscribed sphere of a translational tetrahedron. We note here that these surfaces are important in order to examine and visualize the Dirichlet-Voronoi cells of $\SOL$ geometry. 
It can be assumed by the homogeneity of $\SOL$ that the starting point of a given translation curve segment is $E_0=(1,0,0,0)$.
The other endpoint will be given by its homogeneous coordinates, $P=(1,a,b,c)$. We consider the translation curve segment $t_{E_0P}$ and determine its
parameters $(\phi,\theta,t)$ expressed by the real coordinates $a$, $b$, $c$ of $P$. We obtain directly by equation system (2.12) the following Lemma (see \cite{Sz17-1}):
\begin{lem}[\cite{Sz17-1}]
\begin{enumerate}
\item Let $(1,a,b,c)$ $(b,c \in \bR \setminus \{0\}, a \in \bR)$ be the homogeneous coordinates of the point $P \in \SOL$. The paramerters of the
corresponding translation curve $t_{E_0P}$ are the following
\begin{equation}
\begin{gathered}
\phi=\mathrm{arccot}\Big(-\frac{a}{b} \frac{\mathrm{e}^{c}-1}{\mathrm{e}^{-c}-1}\Big),~\theta=\mathrm{arccot}\Big( \frac{b}{\sin\phi(\mathrm{e}^{c}-1)}\Big),\\
t=\frac{c}{\sin\theta}, ~ \text{where} ~ -\pi < \phi \le \pi, ~ -\pi/2\le \theta \le \pi/2, ~ t\in \bR^+.
\end{gathered} \tag{2.14}
\end{equation}
\item Let $(1,a,0,c)$ $(a,c \in \bR \setminus \{0\})$ be the homogeneous coordinates of the point $P \in \SOL$. The parameters of the
corresponding translation curve $t_{E_0P}$ are the following
\begin{equation}
\begin{gathered}
\phi=0~\text{or}~  \pi, ~\theta=\mathrm{arccot}\Big( \mp \frac{a}{(\mathrm{e}^{-c}-1)}\Big),\\
t=\frac{c}{\sin\theta}, ~ \text{where}  ~ -\pi/2\le \theta \le \pi/2, ~ t\in \bR^+.
\end{gathered} \tag{2.15}
\end{equation}
\item Let $(1,a,b,0)$ $(a,b \in \bR)$ be the homogeneous coordinates of the point $P \in \SOL$. The paramerters of the 
corresponding translation curve $t_{E_0P}$ are the following
\begin{equation}
\begin{gathered}
\phi=\arccos\Big(\frac{x}{\sqrt{a^2+b^2}}\Big),~  \theta=0,\\
t=\sqrt{a^2+b^2}, ~ \text{where}  ~ -\pi < \phi \le \pi, ~ t\in \bR^+.~ ~ \square
\end{gathered} \tag{2.16}
\end{equation}
\end{enumerate}
\end{lem}
This method leads to
\begin{lem}[\cite{Sz17-1}]
The implicit equation of the equidistant surface $\cS_{P_1P_2}(x,y,z)$ of two points $P_1=(1,0,0,0)$ and $P_2=(1,a,b,c)$ in $\SOL$ space:
\begin{enumerate}
\item $c \ne 0$
\begin{equation}\label{bis1}
\begin{gathered}
z\ne 0, c~:~
\frac{|c-z|}{|\mathrm{e}^{c}-\mathrm{e}^z|}\sqrt{(a-x)^2 \mathrm{e}^{2(c+z)}+(\mathrm{e}^{c}-\mathrm{e}^{z})^2+(b-y)^2}=\\
=\frac{|z|}{|\mathrm{e}^{z}-1|}\sqrt{x^2 \mathrm{e}^{2z}+(\mathrm{e}^{z}-1)^2+y^2},\\
z=c~:~\sqrt{(x-a)^2\mathrm{e}^{2c}+(y-b)^2\mathrm{e}^{-2c}}
=\frac{|z|}{|\mathrm{e}^{z}-1|}\sqrt{x^2 \mathrm{e}^{2z}+(\mathrm{e}^{z}-1)^2+y^2},\\
z=0~:~ \frac{|c|}{|\mathrm{e}^{c}-1|}\sqrt{(a-x)^2 \mathrm{e}^{2c}+(\mathrm{e}^{c}-1)^2+(b-y)^2}=\sqrt{x^2+y^2},
\end{gathered} \tag{2.17}
\end{equation}
\item $c=0$
\begin{equation}\label{bis1}
\begin{gathered}
z\ne 0~:~
\frac{|z|}{|\mathrm{e}^z-1|}\sqrt{(a-x)^2 \mathrm{e}^{2z}+(\mathrm{e}^{z}-1)^2+(b-y)^2}=\\
=\frac{|z|}{|\mathrm{e}^{z}-1|}\sqrt{x^2 \mathrm{e}^{2z}+(\mathrm{e}^{z}-1)^2+y^2}\Leftrightarrow \mathrm{e}^{2z}a(a-2x)+b(b-2y)=0,\\
z=0~:~ \sqrt{(x-a)^2+(y-b)^2}=\sqrt{x^2+y^2}\Leftrightarrow xa+yb-\frac{a^2+b^2}{2}.~\square
\end{gathered} \tag{2.18}
\end{equation}
\end{enumerate}
\end{lem}
The process of determining of the center and radius of the circumscribed sphere of a translational tetrahedron 
was examined in the paper \cite{Sz17-1}, here we only recall the necessary results.

We consider $4$ points $A_1$, $A_2$, $A_3$, $A_4$ in the projective model of $\SOL$ space (see Section 2).
These points are the vertices of a {\it translation tetrahedron} in the $\SOL$ space if any two {\it translation segments} connecting the points $A_i$ and $A_j$
$(i<j,~i,j \in \{1,2,3,4\}$) do not have common inner points and any three vertices do not lie in the same translation curve.
Now, the translation segments $A_iA_j$ are called edges of the translation tetrahedron $A_1A_2A_3A_4$.
The circumscribed sphere of a translation tetrahedron is a translation sphere (see Definition 2.2, (2.12) and Fig.~1) that touches each of the tetrahedron's vertices.
As in the Euclidean case the radius
of a translation sphere circumscribed around a tetrahedron $T$ is called the circumradius of $T$, and the center point of this sphere is called the circumcenter of $T$.

\begin{lem}[\cite{Sz17-1}]
For any translation tetrahedron there exists uniquely a translation sphere (called the circumsphere) on which all four vertices lie.
\end{lem}

The procedure to determine the radius and the circumcenter of a given translation tetrahedron is the folowing:

For the circumcenter $C=(1,x,y,z)$ of a given translation tetrahedron $A_1A_2A_3A_4$ $(A_i=(1,x^i,y^i,z^i), ~ i \in \{1,2,3,4\})$
the following system of equations has to hold:
\begin{equation}
d^t(A_1,C)=d^t(A_2,C)=d^t(A_3,C)=d^t(A_4,C). \tag{2.19}
\end{equation}
Therefore it lies on the translation-like bisector surfaces $\cS_{A_i,A_j}$ $(i<j,~i,j \in \{1,2,3,4\}$) which equations are determined in Lemma 2.6.
The coordinates $x,y,z$ of the circumcenter of the circumscribed sphere around the tetrahedron $A_1A_2A_3A_4$ are obtained by the system of equation
derived from these facts:
\begin{equation}
C \in \cS_{A_1A_2}, \cS_{A_1A_3}, \cS_{A_1A_4}. \notag
\end{equation}
Finally, we get the circumradius $r$ as the translation distance e.g. $r=d^t(A_1,C)$.
\begin{figure}[ht]
    \centering
    \includegraphics[width=6cm]{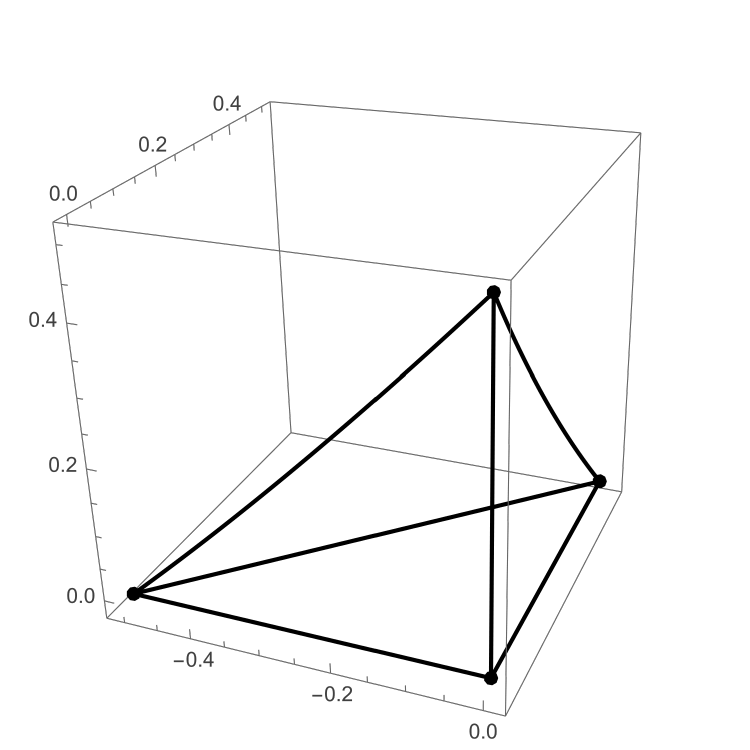} \includegraphics[width=6cm]{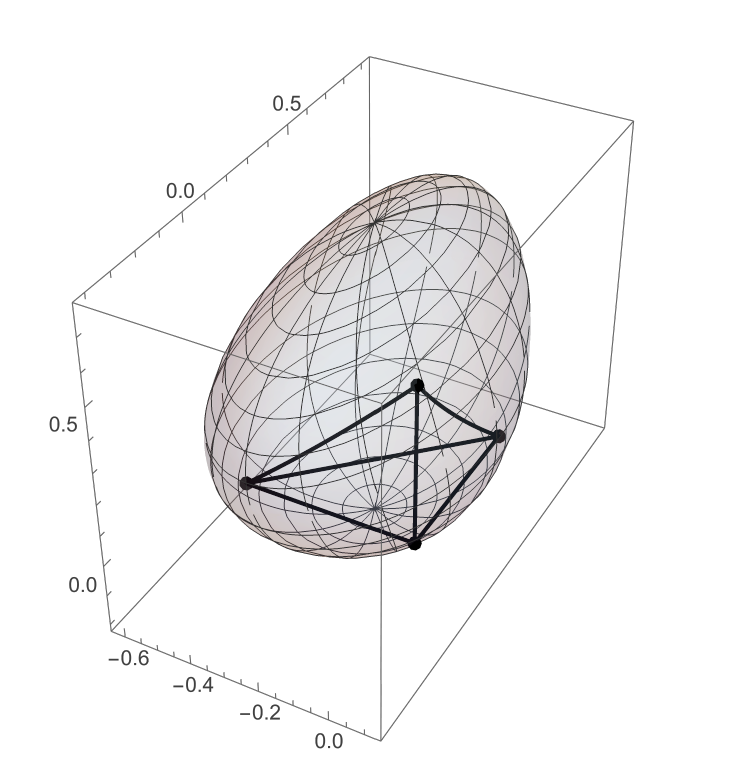} 
    \caption{Tetrahedron with vertices $A_1 = (1,0,0,0), \ A_2 = (1,-0.5,0,0), \
    A_3 = (1,0,0.5,0), \
    A_4 = (1,0,0,0.5)$ and circumradius $r\approx 0.459231$. }
    \label{circumscribedtetrahedron}
\end{figure}
\subsection{On $\SOL$ lattices} 
In this subsection, we recall the necessary results related to $\SOL$ lattices described in the papers \cite{MSz12}, \cite{Sz13-3}, .
\begin{defn}[\cite{MSz12}, \cite{Sz13-3}]
Let $\Gamma < \mathbf{L(R)}$ be a subgroup, generated by three translations 
$\tau_1(t_1^1,t_1^2,t_1^3)$,
$\tau_2(t_2^1,t_2^2,t_2^3)$, $\tau_3(t_3^1,t_3^2,t_3^3)$ with $\mathbf{Z}$ (integer) linear combinations. 
Here upper indices indicate the corresponding $(\mathbf{e}_1,\mathbf{e}_2,\mathbf{e}_3)$ coordinates of basis translations.
$\Gamma$ is called discrete translation group or lattice of $\SOL$, if its action is discrete (by the induced orbit topology), i.e.
there is a compact fundamental "parallelepiped" (with side face 
identifications on its bent side faces) $\widetilde{\cP_\Gamma(N)}=\mathbf{L}/ \Gamma$
(Fig.~2). 
\end{defn}
\begin{thm}[\cite{MSz12}, \cite{Sz13-3}]
Each lattice $\Gamma$ of $\SOL$ has a group presentation (see \cite{S})
\begin{equation}
\begin{gathered}
\Gamma=\Gamma(\Phi)=\{\tau_1,\tau_2,\tau_3: [\tau_1,\tau_2]=1, \tau_3^{-1}\tau_1\tau_3=\tau_1 \Phi^T, \tau_3^{-1}\tau_2\tau_3=\tau_2 \Phi^T
\}, \tag{2.20}
\end{gathered}
\end{equation}
where $t_1^3=t_2^3=0$ holds, and this implies $[\tau_1,\tau_2]=\tau_1^{-1} \tau_2^{-1} \tau_1 \tau_2 = 1$, that means, the commutator of 
$\tau_1, \tau_2$ is the identity, as a first condition of discrete action of $\Gamma$. 
Moreover there exists $\Phi= \begin{pmatrix}
         p & q \\
         r & s \\
\end{pmatrix}\in SL_2(\mathbf{Z})$ with $\mathrm{tr}(\Phi)=p+s > 2, \ ps-qr=1$, such that for above $\tau_1(t_1^1,t_1^2,0), \ \tau_2(t_2^1,t_2^2,0)$
the matrix  $T=\begin{pmatrix}
         t_1^1 & t_1^2 \\
         t_2^1 & t_2^2 \\
\end{pmatrix} \in GL_2(\mathbf{R})$ satisfies the following: 
$T^{-1}\Phi T =:\Phi^T= 
\begin{pmatrix}
         e^{-t_3^3}&0 \\
         0& e^{t_3^3} \\
\end{pmatrix}$ is just a hyperbolic rotation fixed by $t_3^3$ in $\tau_3$ above. Namely,   
$\tau_1 \Phi^T=(t_1^1 e^{-t_3^3}, t_1^2 e^{t_3^3})
\begin{pmatrix}
         \be_1 \\
         \be_2 \\
\end{pmatrix},$
$\tau_2 \Phi^T=(t_2^1 e^{-t_3^3}, t_2^2 e^{t_3^3})
\begin{pmatrix}
         \be_1 \\
         \be_2 \\
\end{pmatrix}$ hold
in the commutative basic vector plane of $\SOL$, spanned by $\be_1$ and $\be_2$. These basis vectors are just the eigenvectors of $\Phi^T$
to eigenvalues $e^{-t_3^3}$ and $e^{t_3^3}$, respectively. 
\end{thm}
\begin{defn}[\cite{MSz12}, \cite{Sz13-3}]
A $\SOL$ point lattice $\Gamma_P(\Phi)$ is a discrete orbit of point $P$ in the $\SOL$ space generated by an arbitrary lattice $\Gamma(\Phi)$ above.
For visualizing a point lattice we have chosen the origin as the starting point, by the homogeneity of $\SOL$.
\end{defn}

By the above considerations and by our modifications, we obtain the following summary:  
\begin{thm}[\cite{MSz12}, \cite{Sz13-3}] Translations $\tau_1, \tau_2, \tau_3$ generate a lattice $\Gamma(\Phi)(\tau_1,\tau_2,\tau_3)$  
in the $\SOL$ space 
if and only if the vectors $\tau_1, \tau_2$ in 
the $[x,y]$ plane generate a lattice, unimodularly equivalent to a $\Gamma(N,p,q)$ lattice; the components $t_3^1,t_3^2 \in \mathbf{R}$ 
are given modulo this 
sublattice $\Gamma^0(\tau_1, \tau_2)$, and  the parameter $t_3^3$ 
in $\tau_3$
satisfies the following (see \cite{MSz12}, \cite{Sz13-3}):
\begin{equation}
     \begin{gathered}
     2 \cosh{t_3^3}=p+s=N, \\
     \begin{pmatrix}
         e^{-t_3^3}&0 \\
         0&e^{t_3^3}
     \end{pmatrix}\Leftrightarrow \begin{pmatrix}
              \cosh{t_3^3}&\sinh{t_3^3} \\
              \sinh{t_3^3}&\cosh{t_3^3}
       \end{pmatrix} \Rightarrow t_3^3=\log(\frac{1}{2}({N+\sqrt{N^2-4})}).  
       \end{gathered} \tag{2.21}
     \end{equation}
\end{thm}
\subsection{On the fundamental parallelepiped}

If we take integers as coefficients, we generate the discrete group $ \langle \tau_1,\tau_2, \tau_3 \rangle $ 
denoted by $\Gamma(\Phi)$, as above. 
We know that the orbit space $\SOL /\Gamma(\Phi)$ is a compact manifold, i.e. a $\SOL$ space form. 

Let $\widetilde{\cP_\Gamma(N)}$ be a {\it fundamental domain} of $\Gamma(\Phi)$ with face identifications. 
The homogeneous coordinates of the vertices of $\widetilde{\cP_\Gamma(N)}$ can be determined in our affine model by the translations
in Definition 2.9 with the parameters $t_i^j, \ i\in\{1,2,3\}, \ j \in \{1,2,3\}$ as follows. (see Fig.~2). 
\begin{equation}
\begin{gathered}
P(1,t_1^1,t_1^2,0), \ P'(1,t_2^1,t_2^2,0), \ P_3(1,t_3^1,t_3^2, t_3^3), \ Q(1,t_1^1+t_2^1,t_1^2+t_2^2,0), \\ 
Q'(1,(t_1^1+t_2^1) e^{-t_3^3}, (t_1^2+t_2^2) e^{t_3^3} ,0), \\
Q^{\tau_3}(1,t_3^1+(t_1^1+t_2^1) e^{-t_3^3}, t_3^2+(t_1^2+t_2^2) e^{t_3^3} ,t_3^3), \
P''(1,t_2^1 e^{-t_3^3}, t_2^2 e^{t_3^3} ,0), \\ P'^{\tau_3}(1,t_3^1+t_2^1 e^{-t_3^3}, t_3^2+t_2^2 e^{t_3^3} ,t_3^3), P^{\tau_3}(1,t_3^1+t_1^1 e^{-t_3^3}, 
t_3^2+t_1^2 e^{t_3^3} ,t_3^3). \tag{2.22}
\end{gathered}
\end{equation}
\begin{defn}[\cite{Sz13-3}]
$\Gamma(\Phi)=\Gamma(\Phi,t_3^1,t_3^2)$ is called a fundamental lattice in $\SOL$ space if $t_3^1=0=t_3^2$ and 
$\{\tau_1,\tau_2=\tau_1 \Phi^T, \tau_3 \}$ 
is a basis of $\Gamma(\Phi)$.
\end{defn}

The following basis describes the {\it type of fundamental lattices}, that means, $\mathrm{tr}(\Phi)=N$ is the characteristic free parameter, and 
\begin{equation}
\frac{t_2^1}{t_1^1}=\frac{N-\sqrt{N^2-4}}{2}=e^{-t_3^3}, \ \ \frac{t_2^2}{t_1^2}=\frac{N+\sqrt{N^2-4}}{2}=e^{t_3^3},  \tag{2.23}
\end{equation}
for $\tau_1(t_1^1,t_1^2,0), \tau_2(t_2^1,t_2^2,0)$.
The case $N=p+s=3$ and $t_1^1=\frac{1}{\sqrt{2}}=t_1^2$, $t_2^1=\frac{3-\sqrt{5}}{2\sqrt{2}},  t_2^2=\frac{3+\sqrt{5}}{2\sqrt{2}}$,
$t_3^1=t_3^2=0, \ t_3^3=\log{\frac{3+\sqrt{5}}{2}}$ is illustrated in Fig.~2.
\begin{figure}[ht]
\centering
\includegraphics[width=12cm]{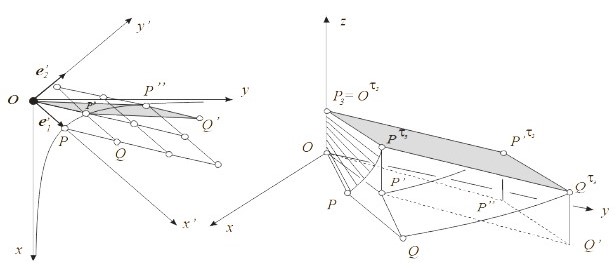} 
\caption{A fundamental parallelepiped for $N=3$.}
\label{}
\end{figure}

It is easy to see by (2.4) and by integration that the volume of $\widetilde{\mathcal{F}}$ can be determined by the following formula:
\begin{thm}[\cite{Sz13-3}]
\begin{equation}
Vol(\widetilde{\mathcal{F}})=\big| \det \begin{pmatrix}
         t_1^1 & t_1^2 \\
         t_1^2 & t_2^2
     \end{pmatrix} \cdot t_3^3 \big| = \big|(t_1^1 t_2^2 -t_1^2 t_2^1)\cdot \log(\frac{1}{2}({N+\sqrt{N^2-4})})\big|. \tag{2.24}
\end{equation}
\end{thm}
\section{The lattice-like translation ball coverings}
To construct a covering of a {\bf fundamental lattice} with congruent translation balls, we want to find a suitable radius, for which placing the centers of the balls at the vertices of the lattice, we can ensure that the balls together cover the whole space. It is also important to try to find a minimal such radius, so that the density of the covering may be minimal.

To ensure that the 8 balls together cover the lattice, we need to make sure that they are convex in an Euclidean sense, which is not true for all spheres in $\SOL$ geometry. 
For this, the following Theorem is used that can be proved using 
classical differential geometry tools, we will not detail this here (see Fig.~3,4,5):

\begin{thm}
    A sphere in the $\SOL$ space is convex in an Euclidean sense, if for its radius: $r \in (0, \frac{\pi}{2}]$. \quad \quad $\square$
\end{thm}

\begin{figure}
    \centering
    \includegraphics[width=9cm]{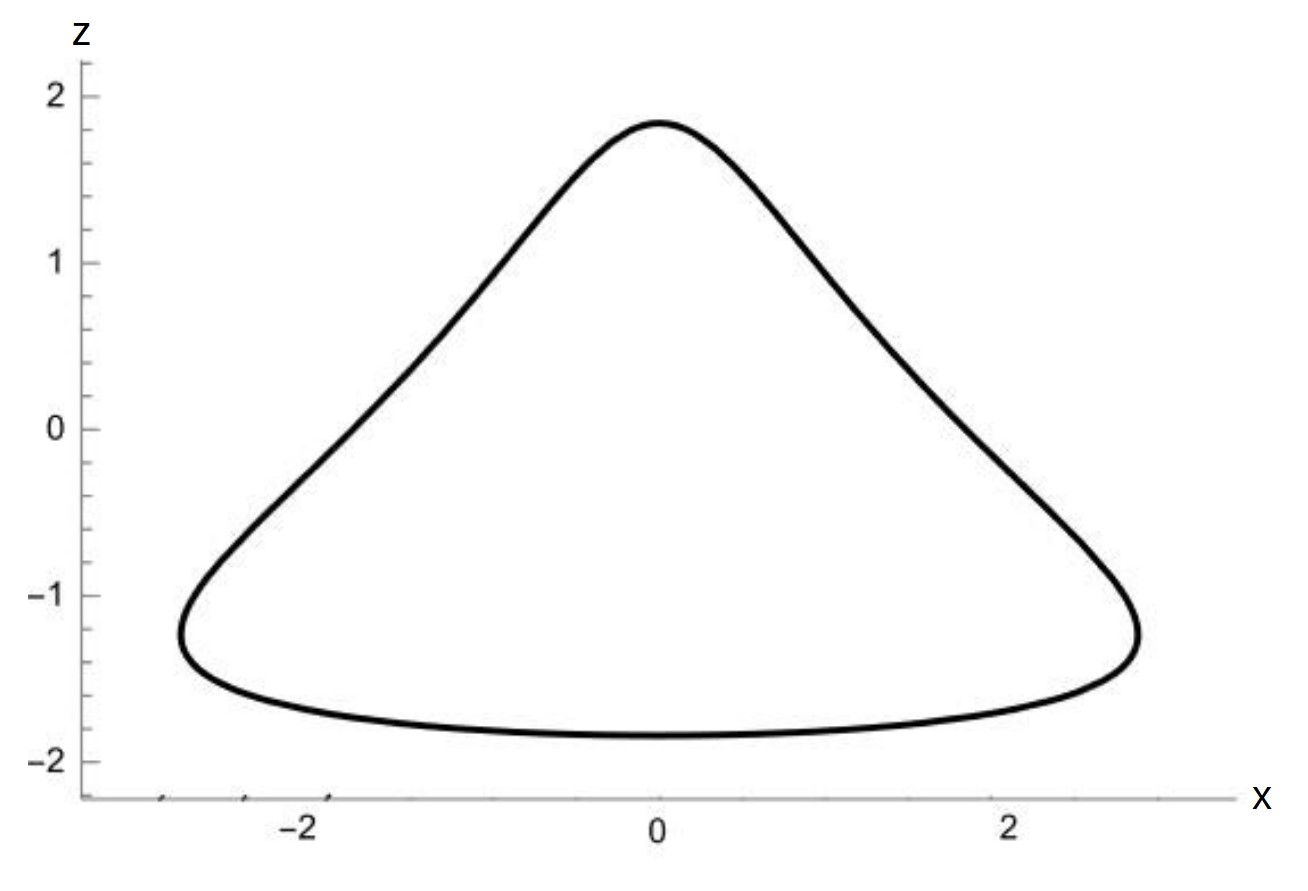}
    \caption{The plane sections of the sphere with $r=2$ parallel to the plane 
    $[y,z]$ at $y=0$}
    \label{gombmetszet}
\end{figure}

\begin{figure}
\centering
        \includegraphics[width=6cm]{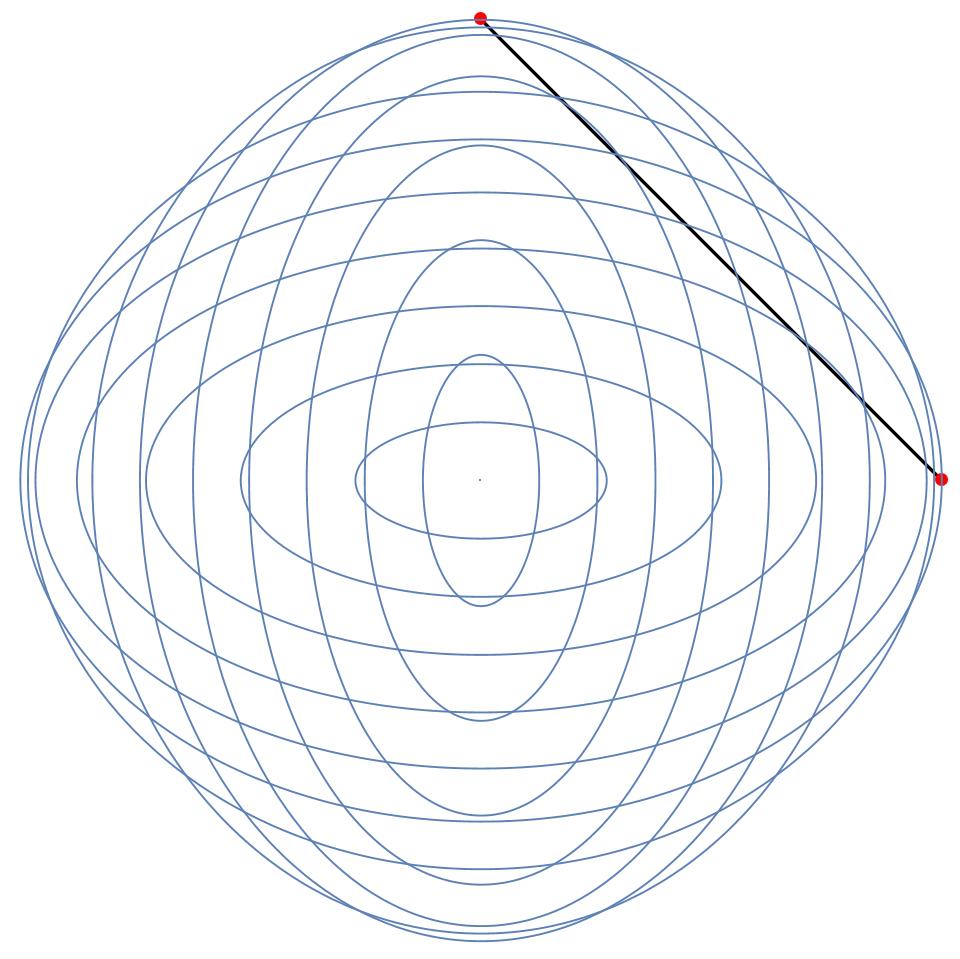} \includegraphics[width=6cm]{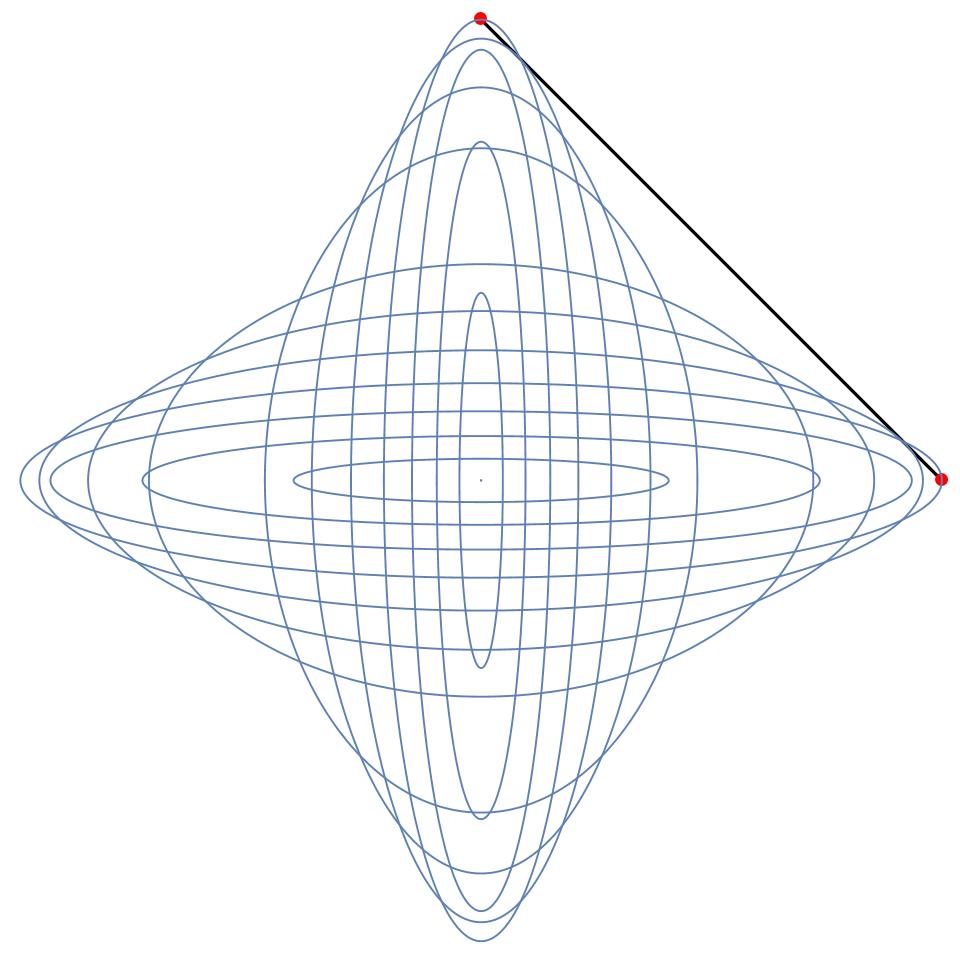}
        
        $r=\frac{\pi}{4}$ \qquad \qquad $r=2.2$ 
        
        \caption{Plane sections of spheres with different radii, parallel to the $[x,y]$ plane.}
    \label{ynemnulla}
\end{figure}

\begin{figure}
\centering
    \includegraphics[width=12cm]{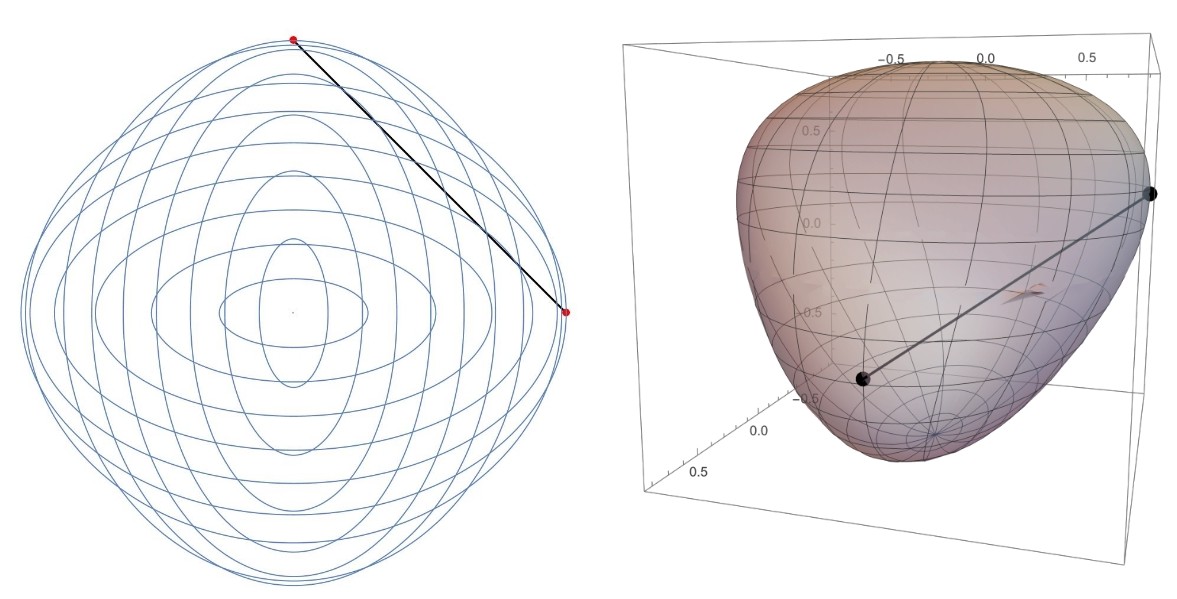} 
    \caption{The visualization of a convex translation ball and 
    its plane sections parallel to the $[x,y]$ base plane.}
    \label{yno}
\end{figure}

\subsection{Coverings of the fundamental lattice}
Let $\cB_\Gamma^C(r)$ be the covering of the $\SOL$ space with translation balls $B^C(r)$ of radius $r$. The centers of the balls form a $\SOL$ point lattice 
$\Gamma(\tau_1,\tau_2,\tau_3)$. We denote an arbitrary $\SOL$ parallelepiped of this lattice 
by $\widetilde{\cP_\Gamma(N)}$, whose 
images by the discrete translation group $\Gamma(\Phi)$ described in subections 
2.3 and 2.4 covers the $\SOL$ space with no gap.

\begin{rem}
    Similarly to the Euclidean space $\bE^d\;(d\geq1)$, an arbitrary 
    lattice $\Gamma$ belongs to a lattice-like covering of equal balls, 
    if the radius of the balls is large enough.
\end{rem}

Taking any translation-like lattice covering $\cB_\Gamma^C(r)$ and shrinking the balls until they don't cover the space any more, 
we get the minimal radius with which the balls still cover the lattice. This defines the least dense covering of the given lattice 
$\Gamma(\tau_1,\tau_2,\tau_3)$ and the threshold value $\Rr$ is called the minimal covering radius.

\begin{gather}
    \Rr := \{r:\;\cB_\Gamma^C(r) \textrm{ is the least dense lattice covering for }\Gamma(\tau_1,\tau_2,\tau_3)\}
    \label{r}
\end{gather}

To determine the the density $\Delta(\cB_\Gamma^C(r))$ of a 
covering, it is sufficient to relate the volume of the minimal covering 
ball and the volume of the solid lattice $\widetilde{\cP_\Gamma(N)}$ (see Theorem 2.13 ), so it can be defined similarly to the Euclidean case:
\begin{defn}
\begin{gather}
    \Delta(\cB_\Gamma^C(r)) := \frac{Vol(\cB_\Gamma^C(r))}{Vol(\widetilde{\cP_\Gamma(N)})}
\end{gather}
\end{defn}
and its minimum is $\Delta(\cB_\Gamma^C(\Rr))$ with radius $R_\Gamma^C$ from formula \ref{r}.

{\it The main problem is finding the lattice $\Gamma(\tau_1,\tau_2, \tau_3)$ to which the optimal minimal density belongs.} 

\begin{equation}
\Delta_{opt}(\mathcal{B}^C_\Gamma)=
\inf_{\Gamma} \Big\{ \Delta(\cB_\Gamma^C(\Rr)) \Big\}. \notag
\end{equation}

\subsection{Thinnest covering}
To determine the thinnest covering for a given fundamental lattice 
$\Gamma(\tau_1,\tau_2,\tau_3)$, we use the following method, 
similar to the one used in \cite{SzV19}.

The lattice is distinctly given by the coordinates $t_1^1, t_1^2$ and the 
parameter $3\leq N\in\NN$. The homogeneous coordinates of the vertices of the 
$\SOL$ parallelepiped 
$\widetilde{\cP_\Gamma(N)}=OPP'QP_3 P^{\tau_3} P'^{\tau_3} Q^{\tau_3}$ can be derived from 
these three attributes (see formula (2.22)). 
It is enough to examine the ball coverings $\cB_\Gamma^C$ of 
$\widetilde{\cP_\Gamma(N)}$.

To find the minimal covering radius $\Rr$ for a given 
lattice $\Gamma(t_1^1, t_1^2, N)$, we first decompose the fundamental 
parallelepiped $\widetilde{\cP_\Gamma(N)}$ into six tetrahedra, which 
together fill it once with no gaps or overlapping: 
$\{O,P_3,P^{\tau_3},P'^{\tau_3}\},\;
	\{O,P,P^{\tau_3},P'\},\\
	\{O,P',P^{\tau_3},P'^{\tau_3}\},\;
	\{Q^{\tau_3},P^{\tau_3},P',P\},\;
	\{Q^{\tau_3},Q,P',P\},\;
	\{Q^{\tau_3},P'^{\tau_3},P',P^{\tau_3}\}$. Then we determine each of 
	their circumscumradii $R_i\; (i=1,2,\dots,6)$ like it is 
	described in subsection 2.2. 
	With the choice $\Rr = \mathrm{max}\{R_i\}$, the lattice-like 
	arrangement of the balls $\cB_\Gamma^C(\Rr)$ cover the fundamental lattice, 
	if the balls are convex in an Euclidean sense, 
	i.e. $\Rr\in(0,\frac{\pi}{2}]$ (see Theorem 3.1).
	
In the table below we have determined the minimum covering radii and the corresponding covering densities for some given lattices and the tetrahedral decomposition given above. 
This obviously may gives different results for different tetrahedral decompositions, but in our experience this decomposition gives the lowest density.
\medbreak
\centerline{\vbox{
\halign{\strut\vrule\quad \hfil $#$ \hfil\quad\vrule
&\quad \hfil $#$ \hfil\quad\vrule &\quad \hfil $#$ \hfil\quad\vrule
\cr
\noalign{\hrule}
\multispan3{\strut\vrule\hfill\bf Some locally optimal lattice-like translation ball coverings. \hfill\vrule}%
\cr
\noalign{\hrule}
\noalign{\vskip5pt}
\noalign{\hrule}
{\mathrm{Lattice ~ parameters}} & \Rr &  \Delta_\Gamma^C  \cr
        \noalign{\hrule}
         t_1^1=0.5,\;t^2_1=0.1,\;N=3 & \approx0.60839 & \approx8.82040\cr
        \noalign{\hrule}
          t_1^1=0.6,\;t^2_1=0.2,\;N=3 & \approx0.76434 & \approx7.31397\cr
        \noalign{\hrule}
          t_1^1=0.7,\;t^2_1=0.2,\;N=3 & \approx0.77643 & \approx6.57346\cr
        \noalign{\hrule}
        t_1^1=0.79,\;t^2_1=0.2,\;N=3 & \approx0.80393 & \approx6.47048\cr
        \noalign{\hrule}
         t_1^1=0.8,\;t^2_1=0.2,\;N=3 & \approx0.80731 & \approx6.47105\cr
        \noalign{\hrule}
        t_1^1=0.9,\;t^2_1=0.2,\;N=3 & \approx0.84427 & \approx6.58552 \cr
        \noalign{\hrule}
        \noalign{\hrule}
         t_1^1=0.5,\;t^2_1=0.1,\;N=4 & \approx0.88821 & \approx13.03880\cr
        \noalign{\hrule}
         t_1^1=0.6,\;t^2_1=0.1,\;N=4 & \approx0.93486 & \approx12.68740\cr
        \noalign{\hrule}
          t_1^1=0.7,\;t^2_1=0.1,\;N=4 & \approx,1.00251 & \approx13.44040\cr
        \noalign{\hrule}
        t_1^1=0.8,\;t^2_1=0.2,\;N=4 & \approx1.34052 & \approx14.24860\cr
        \noalign{\hrule}
        t_1^1=0.9,\;t^2_1=0.2,\;N=4 & \approx1.37453 & \approx13.67550\cr
        \noalign{\hrule}
        \noalign{\hrule}
         t_1^1=0.4,\;t^2_1=0.1,\;N=5 & \approx1.22135 & \approx27.24710\cr
        \noalign{\hrule}
          t_1^1=0.5,\;t^2_1=0.1,\;N=5 & \approx1.20960 & \approx21.16430\cr
        \noalign{\hrule}
          t_1^1=0.6,\;t^2_1=0.1,\;N=5 & \approx 1.22263 & \approx18.22280\cr
        \noalign{\hrule}
        t_1^1=0.7,\;t^2_1=0.1,\;N=5 & \approx 1.32750 & \approx20.08420\cr
        \noalign{\hrule}
         t_1^1=0.8,\;t^2_1=0.1,\;N=5 & \approx1.47274 & \approx24.16320\cr
        \noalign{\hrule}
        }}}
\medbreak

\begin{figure}[h!]
    \centering
    \includegraphics[width=13cm]{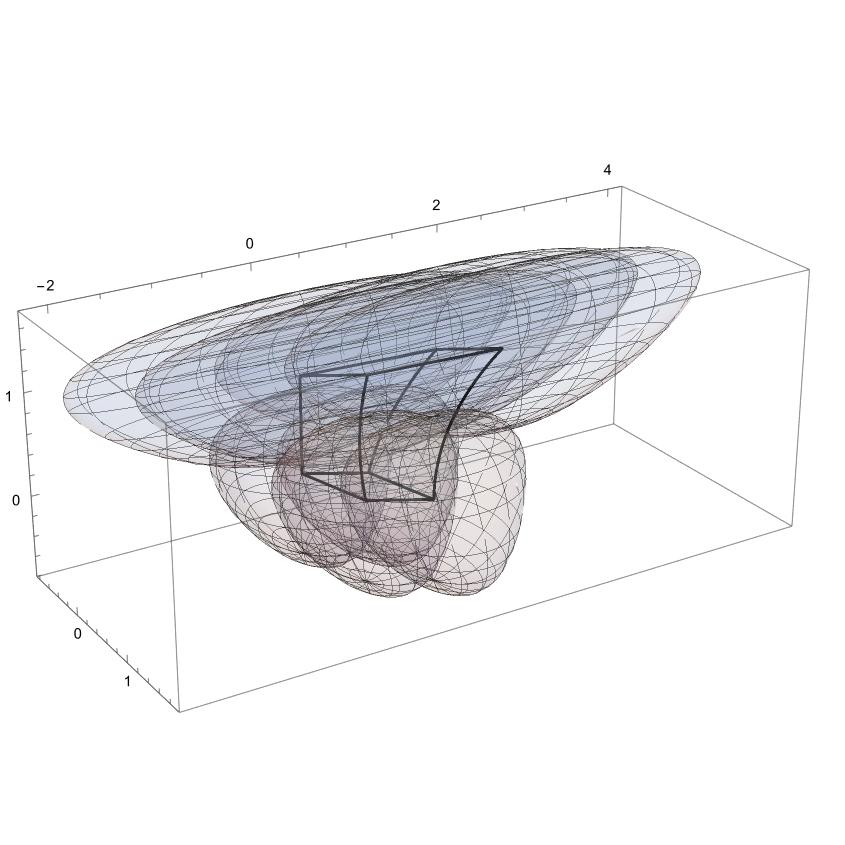}
    \caption{Locally optimal lattice-like translation ball covering of the lattice $\Gamma(0.79,0.2,3)$ with density $\Delta_\Gamma^C\approx6.47048$.}
    \label{locbestcov}
\end{figure}

From the previous computations, by examining many lattices belonging to the fundamental lattices, we obtain the following theorem:
\begin{thm}
    The density of the least dense lattice-like translation ball covering is less than or equal to the locally thinnest covering with congruent 
    translation balls related to the fundamental lattice $\Gamma(t_1^{1,opt}, t_1^{2,opt},N^{opt})$ given by the parameters 
    $t_1^{1,opt}=0.79,\;t^2_{1,opt}=0.2,\;N^{opt}$ (see Fig. \ref{locbestcov}).
    \begin{gather}
        \Delta_{\textrm{opt}}(R_{\textrm{opt}}^C,\tau_1,\tau_2,N) \leq \Delta(R_\Gamma^C,t_1^{1,opt}, t_1^{2,opt},N^{opt}=3)\approx 6,47048.
    \end{gather}
\end{thm}
This covering density can most likely be improved for other lattice types. 
To do this, a similar study should be conducted for the lattices determined in \cite{MSz12}. 
An exact result cannot be expected due to the complexity of the calculations, but approximate better local optima can certainly be achieved. 
These are currently being studied.
\section{The lattice-like cylinder packings and coverings by a type of $\SOL$ translation cylinders}
\subsection{On $\SOL$ cylinders}
\begin{defn}
Let $\cC^i(r)$ be an infinite solid that is bounded by translation curves given by translations parallel to axis $z$ and passing through the points of a 
circle $\cC^b(r)$ of radius $r\in \mathbf{R}^+$ lying 
in the base plane ($[x,y]$ plane, see Section 2) and centred at the origin.
The images $\cC^i_\bt(r)$ of solid $\cC^i(r)$ by $\SOL$ isometries $\bt$ are called {\rm infinite circular cylinders}.
\label{defcyl1}
\end{defn}
The common part of $\cC^i_\bt(r)$ with the base plane is the {\it base figure} of $\cC^i_\bt(r)$ that is denoted by $\cC_\bt(r)$.
\begin{defn}
A {\rm bounded circular cylinder} $\cC(r,h)$ is an isometric image of a solid
which is bounded by the side surface of a infinite circular cylinder $\cC^i(r)$,
its base figure
$\cC^b(r)$ and the translated copy $\cC^c(r)$ of $\cC^b(r)$ by a translation 
parallel to axis $z$ (see formula (2.3)).
The faces $\cC^b(r)$ and $\cC^c(r)$ are called {\rm cover faces}.
The height $h$ (or altitude) of the cylinder is the translation distance between its cover faces.
\label{defcyl2}
\end{defn}
Let us denote the image of the cylinder $\cC^i(r)$ at the $\SOL$ translation $\tau$ 
(see (2.3))
by $\cC^i_\tau(r)$ where the translation is given by parameters $x^0=1,x^1,x^2,x^3=0$ and denote
its common part with the base plane by $\cC_\tau(r)$.

The direct consequence of the (2.2) and (2.3) formulas and the Definitions 4.1-2 is the following;
\begin{lem}
If the translation $\tau$ is given by parameters $x^0=1,x^1,x^2,x^3=0$, then the common part $\cC_\tau(r)$ of $\cC^i_\tau(r)$ with the $[x,y]$ 
plane is a circle in the Euclidean sense. \quad \quad $\square$
\label{lemma"4.3"}
\end{lem} 
Using Lemma 4.3, formula (2.3) and the results of classical differential geometry, we obtain the following;
\begin{rem}
If the base circles $\cC_{\tau_1}^b(r)$ and $\cC_{\tau_2}^b(r)$ of two cylinders touch each other, 
then the corresponding cylinders touch each other along a corresponding translation curve.
\label{cor"4.4"}
\end{rem}
Applying the formulas (2.4), (2.12) we directly obtain  the following lemmas:
\begin{lem}
The perimeter $p(\cC(r))$ and the area ${\mathrm{Area}}(\cC(r))$ of a circle $\cC(r)$ with radius $r$ (Euclidean and the $\SOL$ translation distances are equal in the base plane)
lying in the base plane centred at the origin of the projective model of the $\SOL$
geometry can be calculated with the same formulas in the $\SOL$ and the Euclidean $\EUC$ geometries.
\begin{equation}
p(\cC(r))=2\cdot r \cdot \pi, ~ {\mathrm{Area}}(\cC(r))=r^2 \cdot \pi. \tag{4.1}
\end{equation}
\end{lem}
\begin{cor}
The volume ${\mathrm{Vol}}(\cC(r,h))$ of a bounded fibre-like circular cylinder with radius $r$ and height (or altitude) $h$ in the projective model of the $\SOL$
geometry can be calculated with the following formula
\begin{equation}
{\mathrm{Vol}}(\cC(r,h))=h \cdot r^2 \cdot \pi.\tag{4.2}
\end{equation}
\end{cor}

\begin{figure}[ht]
\centering
       \includegraphics[width=5.5cm]{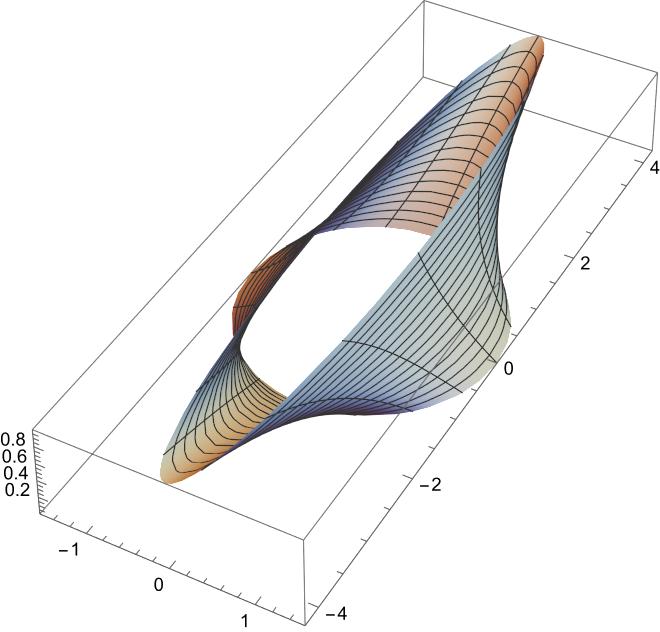} \includegraphics[width=7.5cm]{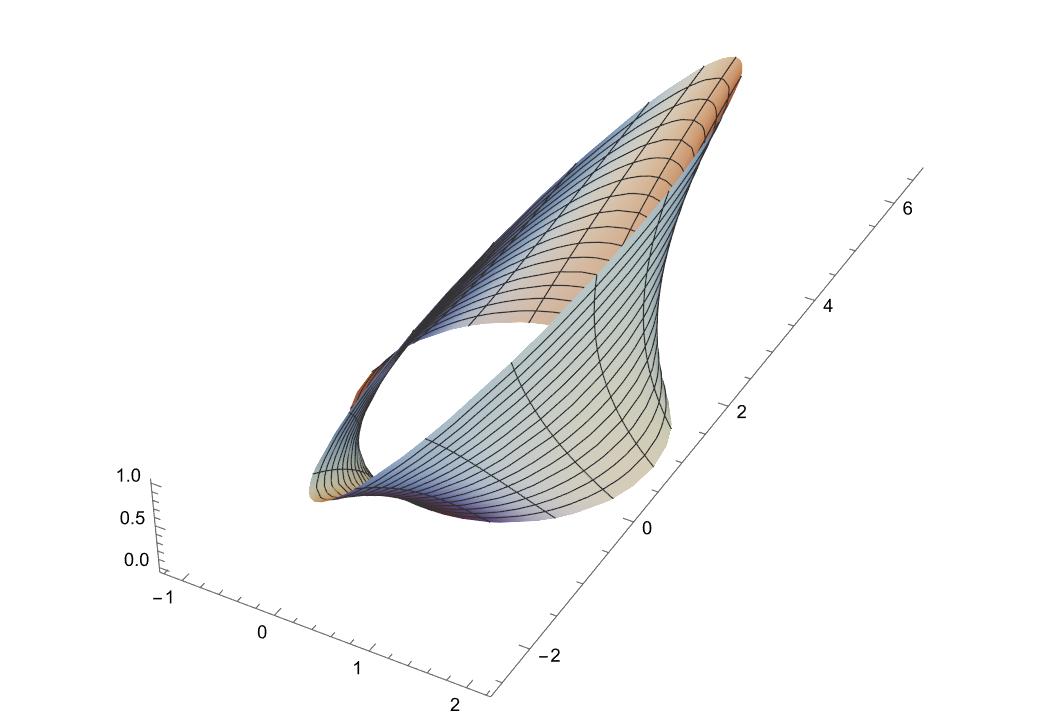} 
       \caption{Translation cylinders with radius $r=\frac{\pi}{2}$, 
       height $h=1$ and different centers in the $\SOL$ space.}
    \label{cyls}
\end{figure}
\subsubsection{Lattice-like cylinder packings, coverings and their densities}
Let us consider a tiling $\cT_\Gamma(N)$ with parallelepipeds related to the fundamental lattices in the $\SOL$ space. Let $\widetilde{\cP_\Gamma(N)}$ 
be one of its tiles with base figure $\widetilde{\cP_\Gamma(N)}^b$.
This tiling is generated by translations $\tau_1, \tau_2, \tau_3$, that determine the lattice $\Gamma(\Phi)(\tau_1,\tau_2,\tau_3)$ too (see Section 2.3-4 and Theorem 2.11).   

{\bf In this section we investigate a large class of lattice-like translation cylinder packings and coverings in $\SOL$ space where the 
tilings are related to the so called {\it fundamental lattices}}.

We consider the corresponding images of the cylinder $\cC^i(r)$ generated by  $\tau_1(t_1^1,t_1^2,0), \ \tau_2(t_2^1,t_2^2,0)$ translations of 
the fundamental lattices (see Theorem 2.11 and formula (2.22)).

By the Lemma 4.3, the common parts $\cC_{\tau_j}^b(r)$ of these images  
with the fundamental plane will be circles in the Euclidean sense. 
Furthermore, the sublattice $\Gamma^0(\tau_1, \tau_2)$ (see Theorems 2.9 and 2.11) of the fundamental lattice in the $[x,y]$ base plane is also 
a lattice in the Euclidean sense. The areas of the circles and subparallelograms also coincide with the Euclidean areas by Theorem 2.13 and Lemma 4.5.

From all this, the following are clearly evident:
\begin{cor}
\begin{enumerate}
\item The infinite circular cylinders $\cC^i_{\tau_j}(r)$ that are determined by the translations of sublattice $\Gamma^0(\tau_1, \tau_2)$ (see Theorem 2.11)
form a lattice-like cylinder packing or covering related to the corresponding fundamental lattice if their common parts with the $[x,y]$ base plane form a 
lattice-like circle packing or covering respectively.
\item If we divide the above infinite cylinders $\cC^i_{\tau_j}(r)$ into bounded circular cylinders with planes parallel to the $[x,y]$ base plane having the height distances 
of the fundamental parallelepiped ($h=t_3^3=\log(\frac{1}{2}({N+\sqrt{N^2-4})})$, see Theorem 2.11) we obtain a lattice like packing or covering with bounded circular cylinders.
\end{enumerate}
\end{cor}

For the density of the optimal lattice like cylinder packing or covering related to the $\cT_\Gamma(N)$ parallelepiped tilings determined by fundamental lattices in the $\SOL$ space 
it is sufficient to relate the volume of the bounded circular cylinder $\cC(r,h)$ ($h=\log(\frac{1}{2}({N+\sqrt{N^2-4})})$)
to that of the parallelepiped $\widetilde{\cP_\Gamma(N)}$. This ratio can be replaced by the ratio of the areas of base figurs ${\mathrm{Area}}(\cC^b(r))$ 
and ${\mathrm{Area}}(\cP_\Gamma(N)^b)$
(see Lemma 4.3 and Lemma 4.5, Corollary 4.4).
\begin{defn}
The density of the lattice-like cylinder packing $\cC\cP(\tau_1,\tau_2,N)$ related to $\cT_\Gamma(N)$ parallelepiped tilings determined by fundamental lattices in 
the $\SOL$ space is the following:
\begin{equation}
\delta(\cC\cP(\tau_1,\tau_2,N)):=\frac{\mathrm{Vol}(\cC(r,h))}{\mathrm{Vol}(\widetilde{\cP_\Gamma(N)})}=\frac{{\mathrm{Area}}(\cC^b(r))}{{\mathrm{Area}}(\widetilde{\cP_\Gamma(N)}^b)} 
\notag
\end{equation}
where $r \in \bR^+$ is the radius of cylinders and $N, \tau_1,\tau_2$ are the corresponding parameters of the fundamental lattice (see Definition 2.12).
\end{defn}
\begin{defn}
The density of the lattice-like cylinder covering $\cC\cC(\tau_1,\tau_2,N)$ related to the $\cT_\Gamma(N)$ parallelepiped tilings determined 
by fundamental lattices in the $\SOL$ space is the following:
\begin{equation}
\Delta(\cC\cC(\tau_1,\tau_2,N)):=\frac{\mathrm{Vol}(\cC(r,h))}{\widetilde{\cP_\Gamma(N)}^b}=\frac{{\mathrm{Area}}(\cC^b(r))}{{\mathrm{Area}}(\widetilde{\cP_\Gamma(N)}^b)} \notag
\end{equation}
where $r \in \bR^+$ is the radius of cylinders and $N, \tau_1,\tau_2$ are 
the corresponding parameters of the fundamental lattice (see Definition 2.12)).
\end{defn}
It follows directly from the above notions and results, that
\begin{cor}
 The densities of cylinder packings cannot exceed the known maximum value of Euclidean circle packing and the optimal densities of 
 cylinder coverings cannot be smaller than the minimum density of 
 Euclidean circle coverings:
 \begin{equation}
 \delta(\cC\cP(\tau_1,\tau_2),N) \leq \frac{\pi}{\sqrt{12}} 
 \approx 0.906900; \quad \Delta(\cC\cC(\tau_1,\tau_2,N)) 
 \ge \frac{2\pi}{\sqrt{27}} \approx 1.209200.\notag
 \end{equation}
\end{cor}
\subsection{Densest packing}
\begin{rem}
    As seen in Subsections 2.3-4, the fundamental lattice 
    $\Gamma(\tau_1,\tau_2,\tau_3)$ is determined by the parameters 
    $\tau_1(t_1^1,t_1^2,0)$ and $N$, so we will be using the notation 
    $\Gamma(t_1^1,t_1^2,N)$ to denote these lattices.
\end{rem}

To determine the densest packing for a given fundamental lattice 
$\Gamma(t_1^1,t_1^2,N)$ we choose the radius with the help of the minimum of the two 
sides and two diagonals of the sublattice 
$\Gamma^0(t_1^1,t_1^2)$; the half of this will be the radius $r_p$ 
of the cylinders. This choice ensures that the base circles, and thus 
the cylinders, do not overlap with each other when placing the centers 
of their base circles to the vertices of the subparallelogram.
\begin{equation}
    r_p = \frac{1}{2}\cdot \min\{|\overline{OP}|,|\overline{OP'}|,|\overline{OQ}|,|\overline{PP'}|\} \quad \text{ (see Fig.~7)}
\end{equation}
To also make sure that the cylinders still do not overlap when 
taking the translates of the subparallelogram, 
we check that $r_p$ is less than or equal to the heights of it. 
If not, then we replace $r_p$ with the smaller of these two heights. 
Clearly, this construction gives optimal packing density for a given lattice.
\medbreak
{\footnotesize{\centerline{\vbox{
\halign{\strut\vrule\quad \hfil $#$ \hfil\quad\vrule
&\quad \hfil $#$ \hfil\quad\vrule &\quad \hfil $#$ \hfil\quad\vrule
\cr
\noalign{\hrule}
\multispan3{\strut\vrule\hfill\bf Some locally optimal lattice-like translation 
cylinder packings. \hfill\vrule}%
\cr
\noalign{\hrule}
\noalign{\vskip5pt}
\noalign{\hrule}
{\mathrm{Lattice ~ parameters}} & r_p &  \delta(\cC\cP(t_1^1,t_1^2,N))  \cr
        \noalign{\hrule}
         t_1^1=1.6,\;t^2_1=1,\;N=3 & \approx 0.943398 & \approx 0.781511\cr
        \noalign{\hrule}
         t_1^1=1.8,\;t^2_1=1.1,\;N=3 & \approx 1.04945 & \approx0.781491\cr
        \noalign{\hrule}
         t_1^1=2.1,\;t^2_1=1.3,\;N=3 & \approx1.23491 & \approx 0.784824 \cr
        \noalign{\hrule}
        t_1^1=2.3,\;t^2_1=1.4,\;N=3 & \approx1.33716 & \approx0.780141 \cr
	 \noalign{\hrule}
	 t_1^1=2.4,\;t^2_1=1.5,\;N=3 & \approx1.41510 & \approx0.781511 \cr
	 \noalign{\hrule}
        t_1^1=2.5,\;t^2_1=1.5,\;N=3 & \approx1.43856 & \approx0.775339 \cr 
        \noalign{\hrule}
        t_1^1=1,\;t^2_1=-\frac{1}{2}+\frac{\sqrt{5}}{2},\; N=3 & \frac{\sqrt{10-2 \sqrt{5}}}{4}, \approx 0.58779 & -\frac{(-5+\sqrt{5}) \pi \sqrt{5}}{(20 \sqrt{5}-20)} \approx 0.78540 \cr 
        \noalign{\hrule}
        }}}}}
\medbreak
From the previous method after careful computations we obtain the following theorem:

\begin{thm}
    The density of the densest lattice-like cylinder packing is greater than or equal to the locally densest packing with congruent 
    translation cylinders related to the fundamental lattice $\Gamma(t_1^{1,opt,p}, t_1^{2,opt,p},N^{opt,p})$ given by the parameters 
    $t_1^{1,opt}=1,\;t^{2,opt}_1=-1/2+(1/2) \sqrt{5}\approx 0.61803,\;N^{opt}=3$ (see Fig.~8).
    \begin{gather}
        0.906900 \approx \frac{\pi}{\sqrt{12}}  \geq
  \delta_{\textrm{opt}}(\cC\cP_{\textrm{opt}}(t_1^1,t_1^2,N))) \geq \notag \\
  \geq \delta(\cC\cP(t_1^{1,opt,p},  
  t_1^{2,opt,p},N^{opt,p}))=-\frac{(-5+\sqrt{5}) \pi \sqrt{5}}{(20 \sqrt{5}-20)}\approx 0.78540. 
    \end{gather}
\end{thm}
\begin{rem}
It is obvious, but can also be seen from the previous table, that the densities do not change when $t_1^1/t_1^2 = c \in \bR^+$ constant.
\end{rem}
\begin{figure}[ht]
    \centering
    \includegraphics[width=13cm]{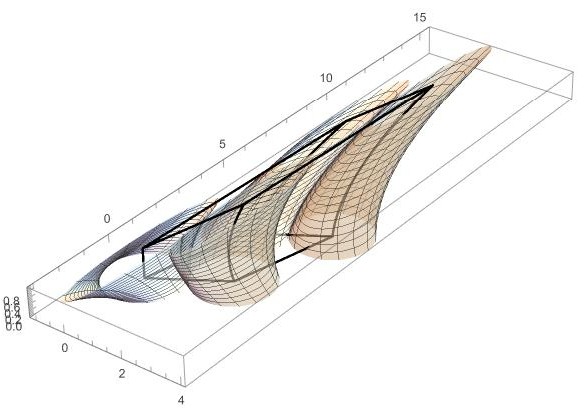}
    \caption{Locally optimal lattice-like cylinder packing of the lattice $\Gamma(t_1^{1,opt},t_1^{2,opt},N^{opt}=3)$ with density $\delta(\cC\cP)\approx0.784824$.}
    \label{locbestpack}
\end{figure}

\subsection{Thinnest covering}
To determine the thinnest covering for a given fundamental lattice 
$\Gamma(t_1^1,t_1^2,N)$, we choose the radius for the cylinder as follows.

We take the sublattice $\Gamma^0(t_1^1,t_1^2)$ and separate it into two triangles $\triangle OPQ, \triangle OQP'$ with one diagonal 
and two other triangles $\triangle OPP', \triangle PP'Q$ with the other diagonal. For these 4 triangles, we calculate the radius of their circumscribed circles. 
We take the maximal radius for both divisions, then the minimum of the two maximums. The radius $r_{c}$ used for the covering will be this, (see Fig.~2 and 9).
\begin{equation}
    r_{c} = \min\{r_{\triangle OPQ},r_{\triangle OPP'}\}
\end{equation}
This choice of $r_c$ provides the locally thinnest cylinder covering for a given lattice. 
\medbreak
{\footnotesize{\centerline{\vbox{
\halign{\strut\vrule\quad \hfil $#$ \hfil\quad\vrule
&\quad \hfil $#$ \hfil\quad\vrule &\quad \hfil $#$ \hfil\quad\vrule
\cr
\noalign{\hrule}
\multispan3{\strut\vrule\hfill\bf Some locally optimal lattice-like translation 
cylinder covering. \hfill\vrule}%
\cr
\noalign{\hrule}
\noalign{\vskip5pt}
\noalign{\hrule}
{\mathrm{Lattice ~ parameters}} & r_c &  \Delta(\cC\cC(t_1^1,t_1^2,N))  \cr
        \noalign{\hrule}
         t_1^1=2.3,\;t^2_1=1,\;N=3 & \approx1.45018 & \approx1.28465   \cr
        \noalign{\hrule}
          t_1^1=2.5,\;t^2_1=1,\;N=3 & \approx1.50151 & \approx1.26701  \cr 
        \hline
          t_1^1=2.6,\;t^2_1=1,\;N=3 & \approx1.52974 & \approx1.26452 \cr
        \noalign{\hrule}
        t_1^1=2.9,\;t^2_1=1.1,\;N=3 & \approx1.69444 & \approx1.26452 \cr
        \hline
         t_1^1=3.1,\;t^2_1=1.2,\;N=3 & \approx1.82991 & \approx 1.26486 \cr
        \noalign{\hrule}
        t_1^1=3.2,\;t^2_1=1.2,\;N=3 & \approx1.85933 & \approx1.26487 \cr
        \noalign{\hrule}
         t_1^1=3.4,\;t^2_1=1.3,\;N=3 & \approx1.99449 & \approx1.26447 \cr 
         \noalign{\hrule}
         t_1^1=1,\;t^2_2=\frac{3-\sqrt{5}}{2},\;N=3 & \frac{(6 \sqrt{5}-15) \sqrt{2}}{-15+5 \sqrt{5}} \approx 0.58632&  -\frac{36}{125}\frac{(2 \sqrt{5}-5)^2 
         \pi \sqrt{5}}{(\sqrt{5}-3)^3} \approx 1.26447 \cr
        \noalign{\hrule}
        }}}}}
\medbreak
From the previous method after careful computations we obtain the following theorem:
\begin{thm}
    The density of the thinnest lattice-like cylinder covering is less than or equal to the locally thinnest covering 
    with congruent translation cylinders related to the fundamental lattice $\Gamma(t_1^{1,opt,c}, t_1^{2,opt,c},N^{opt,c})$ given by the parameters 
    $t_1^{1,opt}=1,\;t^{2,opt}_1=\frac{3-\sqrt{5}}{2} \approx 0.38197,\;N^{opt}=3$ (see Fig.~9).
    \begin{gather}
        1.20920 \approx  \frac{2\pi}{\sqrt{27}} \leq \Delta_{\textrm{opt}}(\cC\cC_{\textrm{opt}}(t_1^1,t_1^2,N)) \leq \notag \\ 
    \leq    \Delta(\cC\cC(t_1^{1,opt,c}, t_1^{2,opt,c},N^{opt,c}))=
        -\frac{36}{125}\frac{(2 \sqrt{5}-5)^2 \pi \sqrt{5}}{(\sqrt{5}-3)^3} \approx 1.26447.
    \end{gather} \quad \quad $\square$
\end{thm}
\begin{figure}[ht]
    \centering
    \includegraphics[width=13cm]{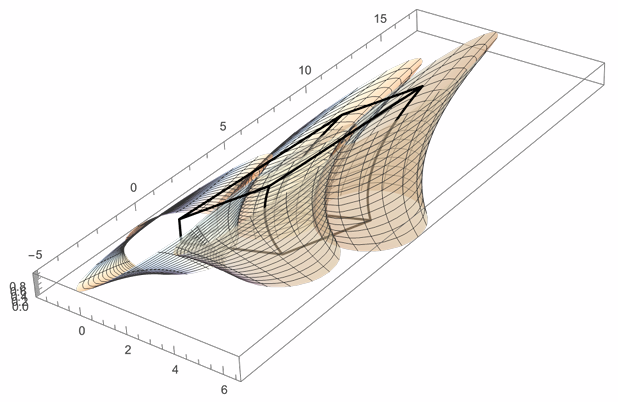}
    \caption{Locally optimal lattice-like cylinder covering of the lattice $\Gamma(t_1^{1,opt,c},t_1^{2,opt,c},N^{opt,c}=3)$ with density $\Delta(\cC\cC) \approx 1.26447$.}
    \label{locbestcovcyl}
\end{figure}
All investigated problems reviewed in this paper can be studied further with the help of other lattices, 
since there are 17 different, nonequivalent types in $\SOL$ (see \cite{MSz12}). 
This may result in finding more optimal, thinner coverings and denser packings for the $\SOL$ space with balls and cylinders.

%

\end{document}